# Classifier Technology and the Illusion of Progress

**David J. Hand**

*Abstract.* A great many tools have been developed for supervised classification, ranging from early methods such as linear discriminant analysis through to modern developments such as neural networks and support vector machines. A large number of comparative studies have been conducted in attempts to establish the relative superiority of these methods. This paper argues that these comparisons often fail to take into account important aspects of real problems, so that the apparent superiority of more sophisticated methods may be something of an illusion. In particular, simple methods typically yield performance almost as good as more sophisticated methods, to the extent that the difference in performance may be swamped by other sources of uncertainty that generally are not considered in the classical supervised classification paradigm.

## 1. INTRODUCTION

In supervised classification, one seeks to construct a rule which will allow one to assign objects to one of a prespecified set of classes based solely on a vector of measurements taken on those objects. Construction of the rule is based on a "design set" or "training set" of objects with known measurement vectors and for which the true class is also known: one essentially tries to extract from the design set the information which is relevant to distinguishing between the classes in terms of the given measurements. It is because the classes are known for the members of this initial data set that the term "supervised" is used: it is as if a "supervisor" has provided these class labels.

Such problems are ubiquitous and, as a consequence, have been tackled in several different research areas, including statistics, machine learning, pattern recognition, computational learning theory and data mining. As a result, a tremendous variety of algorithms and models has been developed for the construction of such rules. A partial list includes linear discriminant analysis, quadratic discriminant analysis, regularized discriminant analysis, the naive Bayes method, logistic discriminant analysis, perceptrons, neural networks, radial basis function methods, vector quantization methods, nearest neighbor and kernel nonparametric methods, tree classifiers such as CART and C4.5, support vector machines and rule-based methods. New methods, new variants on existing methods and new algorithms for existing methods are being developed all the time. In addition, different methods for variable selection, handling missing values and other aspects of data preprocessing multiply the number of tools yet further. General theoretical advances have also been made which have resulted in improved performance at predicting the class of new objects. These include ideas such as bagging, boosting and more general ensemble classifiers. Furthermore, apart from the straightforward development of new rules, theory and practice have been developed for performance assessment. A variety of criteria have been investigated, including measures based on the receiver operating characteristic (ROC) and Brier score, as well as the standard measure of misclassification rate. Subtle estimators of these have been developed, such as jackknife, cross-validation and a variety of bootstrap methods, to overcome the potential optimistic bias which results from simply reclassifying the design set.

*David J. Hand is Professor, Department of Mathematics and Institute for Mathematical Science, Imperial College, Huxley Building, 180 Queen's Gate, London SW7 2AZ, United Kingdom (e-mail: d.j.hand@imperial.ac.uk).*

An examination of recent conference proceedings and journal articles shows that such developments are continuing. In part this is because of new computational developments that permit the exploration of new ideas, and in part it is because of the emergence of new application domains which present new twists on the standard problem. For example, in bioinformatics there are often relatively few cases but many thousands of variables. In such situations the risk of overfitting is substantial and new classes of tools are required. General references to work on supervised classification include [11, 13, 33, 38, 44].

The situation to date thus appears to be one of very substantial theoretical progress, leading to deep theoretical developments and to increased predictive power in practical applications. While all of these things are true, it is the contention of this paper that the practical impact of the developments has been inflated; that although progress has been made, it may well not be as great as has been suggested. The arguments for this assertion are described in the following sections. They develop ideas introduced by Hand [12, 14, 15, 18, 19] and Jamain and Hand [24]. The essence of the argument is that the improvements attributed to the more advanced and recent developments are small, and that aspects of real practical problems often render such small differences irrelevant, or even unreal, so that the gains reported on theoretical grounds, or on empirical comparisons from simulated or even real data sets, do not translate into real advantages in practice. That is, progress is far less than it appears.

These ideas are described in four steps.

First, model-fitting is a sequential process of progressive refinement, which begins by describing the largest and most striking aspects of the data structure, and then turns to progressively smaller aspects (stopping, one hopes, before the process begins to model idiosyncrasies of the observed sample of data rather than aspects of the true underlying distribution). In Section 2 we show that this means that the large gains in predictive accuracy in classification are won using relatively simple models at the start of the process, leaving potential gains which decrease in size as the modeling process is taken further. All of this means that the extra accuracy of the more sophisticated approaches, beyond that attained by simple models, is achieved from "minor" aspects of the distributions and classification problems.

Second, in Section 3 we argue that in many, perhaps most, real classification problems the data points in the design set are not, in fact, randomly drawn from the same distribution as the data points to which the classifier will be applied. There are many reasons for this discrepancy, and some are illustrated. It goes without saying that statements about classifier accuracy based on a false assumption about the identity of the design set distribution and the distribution of future points may well be inaccurate.

Third, when constructing classification rules, various other assumptions and choices are often made which may not be appropriate and which may give misleading impressions of future classifier performance. For example, it is typically assumed that the classes are objectively defined, with no arbitrariness or uncertainty about the class labels, but this is sometimes not the case. Likewise, parameters are often estimated by optimizing criteria which are not relevant to the real aim of classification accuracy. Such issues are described in Section 4 and, once again, it is obvious that these introduce doubts about how the claimed classifier performance will generalize to real problems.

The phenomena with which we are concerned in Sections 3 and 4 are related to the phenomenon of *overfitting*. A model overfits when it models the design sample too closely rather than modeling the distribution from which this sample is drawn. In Sections 3 and 4 we are concerned with situations in which the models may accurately reflect the design distributions (so they do not underfit or overfit), but where they fail to recognize that these distributions, and the apparent classification problems described, are in fact merely a single such problem drawn from a notional distribution of problems. The real aim might be to solve a rather different problem. One might thus describe the issue as one of *problem* uncertainty. To take a familiar example, which we do not explore in detail in this paper because it has been explored elsewhere, the relative *costs* of different kinds of misclassification may differ and may be unknown. A very common resolution is to assume equal costs (Jamain and Hand [24] found that most comparative studies of classification rules made this assumption) and to use straightforward error rate as the performance criterion. However, equality is but one choice, and an arbitrary one at that, and one which we suspect is in fact rarely appropriate. In assuming equal costs, one is adopting a particular problem which may not be the one which is really to be solved. Indeed, things are even worse than this might suggest, because relative misclassification costs may change over time. Provost and Fawcett [36] have described such situations: "Comparison often is difficult in real-world environments because key parameters of the target environment are not known. The optimal cost/benefit tradeoffs and the target class priors seldom are known precisely, and often are subject to change (Zahavi and Levin [47]; Friedman and Wyatt [8]; Klinkenberg and Thorsten [29]). For example, in fraud detection we cannot ignore misclassification costs or the skewed class distribution, nor can we assume that our estimates are precise or static (Fawcett and Provost [6])."

Moving on, our fourth argument is that classification methods are typically evaluated by reporting their performance on a variety of real data sets. However, such empirical comparisons, while superficially attractive, have major problems which are often not acknowledged. In general, we suggest in Section 5 that no method will be universally superior to other methods: relative superiority will depend on the type of data used in the comparisons, the particular data sets used, the performance criterion and a host of other factors. Moreover, the relative performance will depend on the experience the person making the comparison has in using the methods, and this experience may differ between methods: researcher A may find that his favorite method is best, merely because he knows how to squeeze the best performance from this method.

These various arguments together suggest that an apparent superiority in classification accuracy, obtained in "laboratory conditions," may not translate to a superiority in real-world conditions and, in particular, the apparent superiority of highly sophisticated methods may be illusory, with simple methods often being equally effective or even superior in classifying new data points.

## 2. MARGINAL IMPROVEMENTS

This section demonstrates that the extra performance to be achieved by more sophisticated classification rules, beyond that attained by simple methods, is small. It follows that if aspects of the classification problem are not accurately described (e.g., if incorrect distributions have been used, incorrect class definitions have been adopted, inappropriate performance comparison criteria have been applied, etc.), then the reported advantage of the more sophisticated methods may be incorrect. Later sections illustrate how some inaccuracies in the classification problem description can arise.

### 2.1 A Simple Example

Statistical modeling is a sequential process in which one gradually refines the model to provide a better and better fit to the distributions from which the data were drawn. In general, the earlier stages in this process yield greater improvement in model fit than later stages. Furthermore, if one looks at the historical development of classification methods, then the earlier approaches involve relatively simple structures (e.g., the linear forms of linear or logistic discriminant analysis), while more recent approaches involve more complicated structures (e.g., the decision surfaces of neural networks or support vector machines). It follows that the simple approaches will have led to greater improvement in predictive performance than the later approaches which are necessarily trying to improve on the predictive performance obtained by the simpler earlier methods. Put another way, there is a law of diminishing returns.

Although this paper is concerned with supervised classification problems, it is illuminating to examine a simple regression case. Suppose that we have a single response variable $y$ which is to be predicted from $d$ variables $(x_1, \ldots, x_d)^T = \mathbf{x}$. Suppose also that the correlation matrix of $(\mathbf{x}^T, y)^T$ has the form

$$(2.1) \quad \mathbf{\Sigma} = \begin{bmatrix} \Sigma_{11} & \Sigma_{12} \\ \Sigma_{21} & \Sigma_{22} \end{bmatrix} = \begin{bmatrix} (1-\rho)\mathbf{I} + \rho\mathbf{1}\mathbf{1}^T & \boldsymbol{\tau} \\ \boldsymbol{\tau}^T & 1 \end{bmatrix}$$

with $\Sigma_{11} = (1-\rho)\mathbf{I} + \rho\mathbf{1}\mathbf{1}^T$, $\Sigma_{12} = \Sigma_{21}^T = \boldsymbol{\tau}$ and $\Sigma_{22} = 1$, where $\mathbf{I}$ is the $d \times d$ identity matrix, $\mathbf{1} = (1, \ldots, 1)^T$ of length $d$ and $\boldsymbol{\tau} = (\tau, \ldots, \tau)^T$ of length $d$. That is, the correlation between each pair of predictor variables is $\rho$, and the correlation between each predictor variable and the response variable is $\tau$. Suppose also that $\rho, \tau \geq 0$. This condition is not necessary for the argument which follows; it merely allows us to avoid some detail.

Let $V(d)$ be the conditional variance of $y$ given the values of $d$ predictor variables $\mathbf{x}$, as above. Standard results give this conditional variance as

$$(2.2) \qquad V(d) = \Sigma_{22} - \Sigma_{21}\Sigma_{11}^{-1}\Sigma_{12}.$$

Using the result that

$$(2.3) \qquad \begin{aligned} \Sigma_{11}^{-1} &= [(1-\rho)\mathbf{I} + \rho\mathbf{1}\mathbf{1}^T]^{-1} \\ &= \frac{1}{1-\rho}\left\{\mathbf{I} - \frac{\rho\mathbf{1}\mathbf{1}^T}{1+(d-1)\rho}\right\} \end{aligned}$$

[with $-(d-1)^{-1} < \rho < 1$, so that $\Sigma_{11}$ is positive definite], leads to

$$(2.4) \qquad \begin{aligned} V(d) &= 1 - \boldsymbol{\tau}^T \frac{1}{1-\rho}\left\{\mathbf{I} - \frac{\rho\mathbf{1}\mathbf{1}^T}{1+(d-1)\rho}\right\}\boldsymbol{\tau} \\ &= 1 - \frac{d\tau^2}{1-\rho} + \frac{\rho d^2\tau^2}{(1+(d-1)\rho)(1-\rho)}. \end{aligned}$$

From this it follows that the reduction in conditional variance due to adding an extra predictor variable, $x_{d+1}$ (also correlated $\rho$ with the other predictors and $\tau$ with the response variable), is

$$(2.5) \qquad \begin{aligned} X(d+1) &= V(d) - V(d+1) \\ &= \frac{\tau^2}{1-\rho} \\ &\quad + \frac{\rho\tau^2}{1-\rho}\left[\frac{d^2}{1+(d-1)\rho} - \frac{(d+1)^2}{1+d\rho}\right]. \end{aligned}$$

Note that the condition $-(d-1)^{-1} < \rho < 1$ must still be satisfied when $d$ is increased.

Now consider two cases:

*Case* 1. When the predictor variables are uncorrelated, $\rho = 0$. From (2.5), we obtain $X(d+1) = \tau^2$. That is, if the predictor variables are mutually uncorrelated and each has correlation $\tau$ with the response variable, then each additional predictor reduces the variance of the conditional variance of $y$ given the predictors by $\tau^2$. [Of course, by setting $\rho = 0$ in (2.4) we see that this is only possible up to $d = \tau^{-2}$ predictors. With this many predictors the conditional variance of $y$ given $\mathbf{x}$ has been reduced to zero.]

*Case* 2. $\rho > 0$. Plots of $V(d)$ for $\tau = 0.5$ and for a range of $\rho$ values are shown in Figure 1. When there is reasonably strong mutual correlation between the predictor variables, the earliest ones contribute substantially more to the reduction in variance remaining unexplained than do the later ones. The case $\rho = 0$ consists of a diagonal straight line running from 1 down to zero. In the case $\rho = 0.9$, almost all of the variance in the response variable is explained by the first chosen predictor.

This example shows that the reduction in conditional variance of the response variable decreases with each additional predictor we add, even though each predictor has an identical correlation with the response variable (provided this correlation is greater than 0). The reason for the reduction is, of course, the mutual correlation between the predictors: much of the predictive power of a new predictor has already been accounted for by the existing predictors.

In real applications, the situation is generally even more pronounced than in this illustration. Usually, in real applications, the predictor variables are not identically correlated with the response, and the predictors are selected sequentially, beginning with those which maximally reduce the conditional variance. In a sense, then, the example above provides a lower bound on the phenomenon: in real applications the proportion of the gains attributable to the early steps is even greater.

### 2.2 Decreasing Bounds on Possible Improvement

We now return to supervised classification. For illustrative purposes, suppose that misclassification rate is the performance criterion, although similar arguments apply with other criteria. Ignoring issues of overfitting, adding additional predictor variables can only lead to a decrease in misclassification rate.

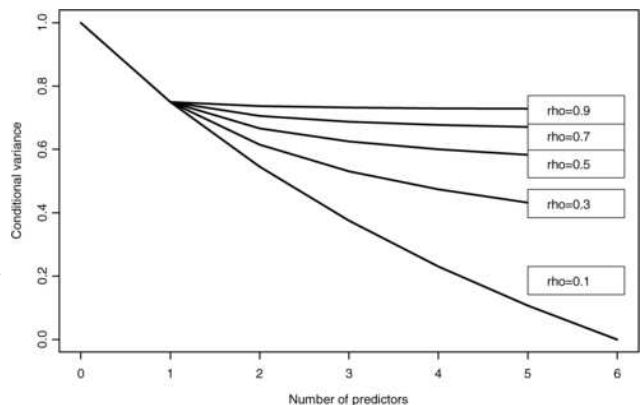

FIG. 1. *Conditional variance of response variable as additional predictors are added for $\tau = 0.5$. A range of values of $\rho$ is shown.*

The simplest model is that which uses no predictors, leading, in the two-class case, to a misclassification rate of $m_0 = \pi_0$, where $\pi_0$ is the prior probability of the smaller class. Suppose that a predictor variable is now introduced which has the effect of reducing the misclassification rate to $m_1 < m_0$. Then the scope for further improvement is only $m_1$, which is less than the original scope $m_0$. Furthermore, if $m_1 < m_0 - m_1$, then all future additions necessarily improve things by less than the first predictor variable. In fact, things are even more extreme than this: one cannot further reduce the misclassification rate by more than $m_1 - m_b$, where $m_b$ is the Bayes error rate. To put it another way, at each step the maximum possible increase in predictive power decreases, so it is not surprising that, in general, at each step the additional contribution to predictive power decreases.

### 2.3 Effectiveness of Simple Classifiers

Although the literature contains examples of artificial data which simple models cannot separate (e.g., intertwined spirals or checkerboard patterns), such data sets are exceedingly rare in real life. Conversely, in the two-class case, although few real data sets have exactly linear decision surfaces, it is common to find that the centroids of the predictor variable distributions of the classes are different, so that a simple linear surface can do surprisingly well as an estimate of the true decision surface. This may not be the same as "can do surprisingly well in classifying the points," since in many problems the Bayes error rate is high, meaning that no decision surface can separate the distributions of such problems very well. However, it means that the dramatic steps in improvement in classifier accuracy are made in the simple first steps. This is a phenomenon which has been noticed by others (e.g., Rendell and Seshu [37]; Shavlik, Mooney and Towell [41]; Mingers [34]; Weiss, Galen and Tadepalli [45]; Holte [22]). Holte [22], in particular, carried out an investigation of this phenomenon. His "simple classifier" (called $1R$) consists of a partition of a single variable, with each cell of the partition possibly being assigned to a different class: it is a multiple-split single-level tree classifier. A search through the variables is used to find that which yields the best predictive accuracy. Holte compared this simple rule with C4.5, a more sophisticated tree algorithm, finding that "on most of the datasets studied, $1R$'s accuracy is about 3 percentage points lower than C4's."

We carried out a similar analysis. Perhaps the earliest classification method formally developed is Fisher's linear discriminant analysis [7]. Table 1 shows misclassification rates for this method and for the best performing method we could find in a search of the literature (these data were abstracted from the data accumulated by Jamain [23] and Jamain and Hand [24]) for a randomly selected sample of ten data sets. The first numerical column shows the misclassification rate of the best method we found $(m_T)$, the second shows that of linear discriminant analysis $(m_L)$, the third shows the default rule of assigning every point to the majority class $(m_0)$ and the final column shows the proportion of the difference between the default rule and the best rule which is achieved by linear discriminant analysis $[(m_0 - m_L)/(m_0 - m_T)]$. It is likely that the best rules, being the best of rules which many researchers have applied, are producing results near the Bayes error rate.

TABLE 1
*Performance of linear discriminant analysis and the best result we found on ten randomly selected data sets*

| Data set | Best method e.r. | Lindisc e.r. | Default rule | Prop linear |
|---|---|---|---|---|
| Segmentation | 0.0140 | 0.083 | 0.760 | 0.907 |
| Pima | 0.1979 | 0.221 | 0.350 | 0.848 |
| House-votes16 | 0.0270 | 0.046 | 0.386 | 0.948 |
| Vehicle | 0.1450 | 0.216 | 0.750 | 0.883 |
| Satimage | 0.0850 | 0.160 | 0.758 | 0.889 |
| Heart Cleveland | 0.1410 | 0.141 | 0.560 | 1.000 |
| Splice | 0.0330 | 0.057 | 0.475 | 0.945 |
| Waveform21 | 0.0035 | 0.004 | 0.667 | 0.999 |
| Led7 | 0.2650 | 0.265 | 0.900 | 1.000 |
| Breast Wisconsin | 0.0260 | 0.038 | 0.345 | 0.963 |

The striking thing about Table 1 is the large values of the percentages of classification accuracy gained by simple linear discriminant analysis. The lowest percentage is 85% and in most cases over 90% of the achievable improvement in predictive accuracy, over the simple baseline model, is achieved by the simple linear classifier.

I am grateful to Willi Sauerbrei for pointing out that when the error rates of both the best method and the linear method are small, the large proportion in achievable accuracy which can be obtained by the linear method corresponds to the error rate of the linear method being a large multiple of that of the best method. For example, in the most extreme case in Table 1, the results for the segmentation data show that the linear discrimination error rate is nearly six times that of the best method. On the other hand, when the error rates are small, this large difference will correspond to only a small proportion of new data points. Small differences in error rate are susceptible to the issues raised in Sections 3 and 4: they may vanish when problem uncertainties are taken into account.

## 2.4 The Flat Maximum Effect

Even within the context of classifiers defined in terms of simple linear combinations of the predictor variables, it has often been observed that the major gains are made by (for example) weighting the variables equally, with only little further gains to be had by careful optimization of the weights. This phenomenon has been termed the *flat maximum* effect [13, 43]: in general, often quite large deviations from the optimal set of weights will yield predictive performance not substantially worse than the optimal weights. An informal argument that shows why this is often the case is as follows.

Let the predictor variables be $(x_1,\ldots,x_d)^T = \mathbf{x}$ and, for simplicity, assume that $E(x_i) = 0$ and $V(x_i) = 1$ for $i = 1,\ldots,d$. Let $\mathbf{\Sigma} = \{r\}_{ij}$ be the correlation matrix between these variables. Now define two weighted sums

$$w = \sum_{i=1}^{d} w_i x_i \quad \text{and} \quad v = \sum_{i=1}^{d} v_i x_i,$$

using respective weight vectors $(w_1,\ldots,w_d)$ and $(v_1,\ldots,v_d)$. In general, $r(w,v)$, the correlation between $w$ and $v$, can take extreme values of $+1$ and $-1$, but suppose we restrict the weights to be nonnegative, $w_i, v_i \geq 0$ for $i = 1,\ldots,d$, and also require $\sum w_i = 1$ and $\sum v_i = 1$. Using these conditions, a little algebra shows that

$$r(v,w) \geq \sum_i \sum_j v_i w_j r(x_i, x_j).$$

Now, with equal weights, $v_i = 1/d, i = 1,\ldots,n$, we obtain

$$r(v,w) \geq \frac{1}{d} \sum_i \sum_j w_j r(x_i, x_j)$$

$$\geq \frac{1}{d} \sum_i \sum_j w_j r(x_i, x_k),$$

where $k = \arg\min_j r(x_i, x_j)$.

From this,

$$r(v,w) \geq \frac{1}{d} \sum_i \sum_j w_j r(x_i, x_k)$$

$$= \frac{1}{d} \sum_i r(x_i, x_k).$$

In words, the correlation between an arbitrary weighted sum of the $x$ variables (with weights summing to 1) and the simple combination using equal weights is bounded below by the smallest row average of the entries in the correlation matrix of the $x$ variables. Hence if the correlations are all high, the simple average will be highly correlated with any other weighted sum: the choice of weights will make little difference to the scores. The gain to be made by the extra effort of optimizing the weights may not be worth the effort.

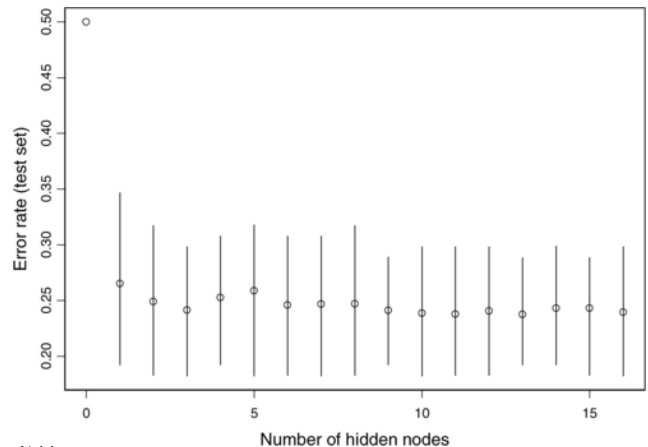

FIG. 2. *Effect on misclassification rate of increasing the number of hidden nodes in a neural network to predict the class of the sonar data.*

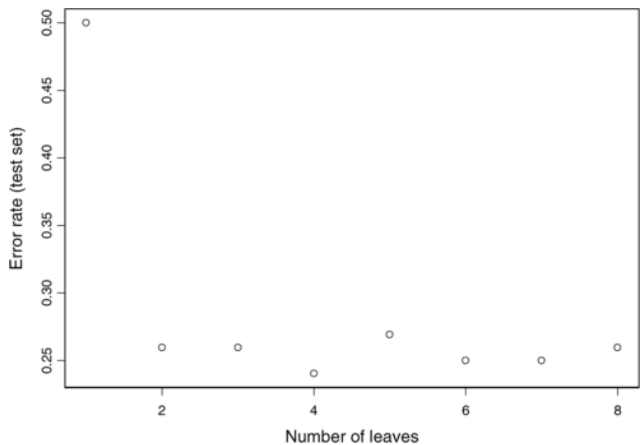

Fig. 3. *Effect on misclassification rate of increasing the number of leaves in a tree classifier to predict the class of the sonar data.*

## 2.5 An Example

As a simple illustration of how increasing model complexity leads to a decreasing rate of improvement, we fitted models to the sonar data from the University of California, Irvine (UCI) data base. This data set consists of 208 observations, 111 of which belong to the class "metal" and 97 of which belong to the class "rock." There are 60 predictor variables. The data were randomly divided into two parts, and a succession of neural networks with increasing numbers of hidden nodes was fitted to half of the data, with the other half being used as a test set. The error rates are shown in Figure 2. The left-hand point, corresponding to 0 nodes, is the baseline misclassification rate achieved by assigning everyone in the test set to the larger class. The error bars are 95% confidence intervals calculated from 100 networks in each case. Figure 3 shows a similar plot, but this time for a recursive partitioning tree classifier applied to the same data. The horizontal axis shows increasing numbers of leaf nodes. Standard methods of tree construction were used, in which a large tree is pruned back to the requisite number of nodes. In both of these figures we see the dramatic improvement arising from fitting the first nontrivial model. This far exceeds the subsequent improvement obtained in any later step.

## 3. DESIGN SAMPLE SELECTION

Intrinsic to the classical supervised classification paradigm is the assumption that the data in the design set are randomly drawn from the same distribution as the points to be classified in the future. Sometimes slight variants of the sampling scheme are used, for example, drawing samples separately from each class, but the assumption that future points to be classified are drawn from the same distributions as the design set is always made. Unfortunately, as we illustrate in this section, there are several reasons why this assumption may not be justified. In fact, as with our suggestion that the common choice of equal misclassification costs may be more often inappropriate than appropriate, we suspect that the assumption that the design distribution is representative of the distribution from which future points will be drawn is perhaps more often incorrect than correct.

If the distribution underlying the design data and that underlying future points to be classified do differ, then elaborate optimization of the classifier using the design data may be wasted effort: the performance difference between two classifiers may be irrelevant in the context of the differences arising between the design and future distributions. In particular, we suggest, more sophisticated classifiers, which almost by definition model small idiosyncrasies of the distribution underlying the design set, will be more susceptible to wasting effort in this way: the grosser features of the distributions (modeled by simpler methods) are more likely to persist than the smaller features (modeled by the more elaborate methods).

## 3.1 Population Drift

A fundamental assumption of the classical paradigm is that the various distributions involved do not change over time. In fact, in many applications this is unrealistic and the population distributions are nonsta-tionary. For example, it is unrealistic in most commercial applications concerned with human behavior: customers will change their behavior with price changes, with changes to products, with changing competition and with changing economic conditions. Hoadley [21] remarked "the test sample is supposed to represent the population to be encountered in the future. But in reality, it is usually a random sample of the current population. High performance on the test sample does not guarantee high performance on future samples, things do change" and "there is always a chance that a variable and its relationships will change in the future. After that, you still want the model to work. So don't make any variable dominant." He is cautioning against making the model fit the design distribution too well. The last point about not making any variable dominant is related to the flat maximum effect, described above.

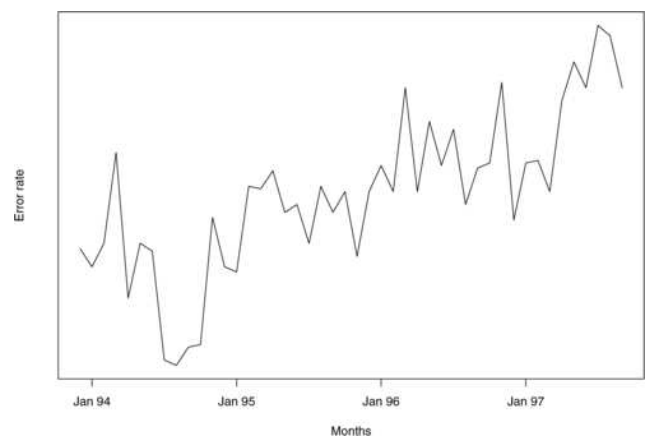

Fig. 4. *Evolution of misclassification rate of a classifier built at the start of the period.*

Among the most important reasons for changes to the distribution of applicants are changes in marketing and advertising practices. Changes to the distributions that describe the customers explain why, in the credit scoring and banking industries [16, 20, 39, 42], the classification rules used to predict which applicants are likely to default on loans are updated every few months: their performance degrades, not because the rules themselves change, but because the distributions to which they are being applied change [27].

An example of this is given in Figure 4. The available data consisted of the true classes ("bad" or "good") and the values of 17 predictor variables for 92,258 customers taking out unsecured personal loans with a 24-month term given by a major UK bank during the period 1 January 1993 to 30 November 1997; 8.86% of the customers belonged to the bad class. The figure shows how the misclassification rate for a classification rule built on data just preceding the start of the displayed period changed over time. Since the coefficients of the classifier were not changing, the deterioration in performance must be due to shifts in the distributions of customers over time.

An illustration of how this "population drift" phenomenon affects different classifiers differentially is given in Figure 5. For the purposes of this illustration we used a linear discriminant analysis (LDA) as a simple classifier and a tree model as a more complicated classifier. For the design set we used customers $1, 3, 5, 7, \ldots, 4999$. We then applied the classifiers to alternate customers, beginning with the second, up to the 60,000th customer. This meant that different customers were used for designing and testing, even during the initial period, so that there would be no overfitting in the reported results. Figure 5 shows lowess smooths of the misclassification cost [i.e., misclassification rate, with customers from each class weighted so that $c_0/c_1 = \pi_1/\pi_0$, where $c_i$ is the cost of misclassifying a customer from class $i$ and $\pi_i$ is the prior (class size) of class $i$]. As can be seen from the figure, the tree classifier (the lower curve) is initially superior (has smaller loss), but after a time its superiority begins to fade. Superficial examination of the figure might suggest that the effect takes a long time to become apparent, not really manifesting itself until around the 40,000th customer, but consider that, in an application such as this, *the data are always retrospective*. In the present case, one cannot determine the true class until the entire 24 month loan term has elapsed. [In fact, of course, this is not quite true: if a customer defaults before the end of the term, then their class (bad) is known, but otherwise their true (good or bad) class is not known until the end, so that to obtain an unbiased sample, one has to wait until the end. Survival analysis models can be constructed to allow for this, but that is leading us away from the point.] For our problem, to accumulate an unbiased sample of 5000 customers with known true outcome, one would have to wait until two years after the 5000th customer had been accepted. In terms of the horizontal axis in Figure 5, this means that the model would be built, and would be initially used at around the time that the 40,000th customer was being considered. The figure shows that this is just when the model degrades. The changes in population structure which occurred during the two years which elapsed while we waited for the true classes of the 5000 design set customers to become known have reduced any advantage that the more sophisticated tree model may have.

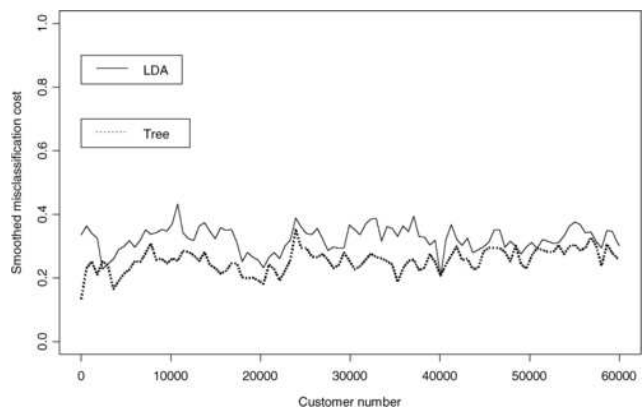

Fig. 5. *Lowess smooths of cost-weighted misclassification rate for a tree model and LDA applied to customers* $2, 4, 6, \ldots, 60{,}000$.

In summary, the apparent superiority of the more sophisticated tree classifier over the very simple linear discriminant classifier is seen to fade when we take into account the fact that the classifiers must necessarily be applied in the future to distributions which are likely to have changed from those which produced the design set. Since, as demonstrated in Section 2, the simple linear classifier captures most of the separation between the classes, the additional distributional subtleties captured by the tree method become less and less relevant when the distributions drift. Only the major aspects are still likely to hold.

The impact of population drift on supervised classification rules is nicely described by the American philosopher Eric Hoffer, who said, "In times of change, learners inherit the Earth, while the learned find themselves beautifully equipped to deal with a world that no longer exists."

### 3.2 Sample Selectivity Bias

The previous subsection considered the impact on classification rules of distributions which changed over time. There is little point in optimizing the rule to the extent that it models aspects of the distributions and decision surface which are likely to have changed by the time the rule is applied. Similar futility applies if a selection process means that the design sample is drawn from a distribution distorted in some way from that to which the classification rule is to be applied. In fact, I suspect that this may be common. Consider, for example, a classification rule aimed at differential medical diagnosis or medical screening. The rule will have been developed on a sample of cases (including members of each class). Perhaps these cases will be drawn from a particular hospital, clinic or health district. Now all sorts of demographic, social, economic and other factors influence who seeks and is accepted for treatment, how severe the cases being treated are, how old they are and so on. In general, it would be risky to assume that these selection criteria are the same for all hospitals, clinics or health districts. This means that the fine points of the classification rule are unlikely to hold. One might expect its coarser features to be true across different such sets of cases, but the detailed aspects will reflect particular properties of the population from which the design data were drawn. In fact, there are some subtleties here. Suppose that the classification rule follows the diagnostic paradigm [directly modeling $p(c|\mathbf{x})$, the probability of class membership, $c$, given the descriptor vector $\mathbf{x}$], rather than the sampling paradigm [which models $p(c|\mathbf{x})$ indirectly from the $p(\mathbf{x}|c)$ using Bayes' theorem]. Then if $\mathbf{x}$ spans the space of all predictors of class membership and if the model form chosen for $p(c|\mathbf{x})$ includes the "true" model, then sampling distortions based on $\mathbf{x}$ alone will not adversely influence the classifier: the classifier built in one clinic will also apply elsewhere. Of course, it would be a brave person who could confidently assert that these two conditions held. Such subtleties aside, what this means, again, is that effort spent on overrefining the classification model is probably wasted effort and, in particular, that fine differences between different classification rules should not be regarded as carrying much weight.

This problem of sample selection and how it might be tackled has been the subject of intensive research, especially by the medical statistics and econometrics communities, but appears not to have been of great concern to researchers on classification methods. Having said that, one area that involves sample selectivity in classification problems which has attracted research interest arises in the retail financial services industry, as in the previous section. Here, as in that section, the aim is to predict, for example, on the basis of application and other background variables, whether or not an applicant is likely to be a good customer. Those expected to be good are accepted, and those expected to be bad are rejected. For those that have been accepted, we subsequently discover their true good or bad class. For the rejected applicants, however, we never know whether they are good or bad. The consequence is that the resulting sample is distorted as a sample from the population of applicants, which is our real interest for the future. Measuring the performance or attempting to build an improved classification rule using those individuals for which we do know the true class (which is needed for supervised classification) has the potential to be highly misleading for the overall applicant population. In particular, it means that using highly sophisticated methods to squeeze subtle information from the design data is pointless. This problem is so ubiquitous in the personal financial services sector that it has been given its own name—*reject inference* [17].

## 4. PROBLEM UNCERTAINTY

Section 3 looked at mismatches between the distributions modeled by the classification rule and the distributions to which it was applied. This is an obvious way in which things may go awry, but there are many others, perhaps not so obvious. This section illustrates just three.

### 4.1 Errors in Class Labels

The classical supervised classification paradigm is based on the assumption that there are no errors in the true class labels. If one expects errors in the class labels, then one can attempt to build models which explicitly allow for this, and there has been work to develop such models. Difficulties arise, however, when one does not expect such errors, but they nevertheless occur.

Suppose that, with two classes, the true posterior class probabilities are $p(1|\mathbf{x})$ and $p(2|\mathbf{x})$, and that a (small) proportion $\delta$ of each class is incorrectly believed to come from the other class at each $\mathbf{x}$. Denoting the apparent posterior probability of class 1 by $p^*(1|\mathbf{x})$, we have

$$p^*(1|\mathbf{x}) = (1-\delta)p(1|\mathbf{x}) + \delta p(2|\mathbf{x}).$$

It follows that if we let $r(\mathbf{x}) = p(1|\mathbf{x})/p(2|\mathbf{x})$ denote the true odds and let $r^*(\mathbf{x}) = p^*(1|\mathbf{x})/p^*(2|\mathbf{x})$ denote the apparent odds, then

$$(4.1) \qquad r^*(\mathbf{x}) = \frac{r(\mathbf{x}) + \varepsilon}{\varepsilon r(\mathbf{x}) + 1}$$

with $\varepsilon = \delta/(1-\delta)$.

With small $\varepsilon$, (4.1) is monotonic increasing in $r(\mathbf{x})$, so that contours of $r(\mathbf{x})$ map to corresponding contours of $r^*(\mathbf{x})$. In particular, if the true optimal decision surface is $r(\mathbf{x}) = k$ ($k$ is determined by the relative misclassification costs), then the optimal decision surface when errors are present is given by $r^*(\mathbf{x}) = k^*$, with $k^* = (k+\varepsilon)/(\varepsilon k + 1)$. Unfortunately, if the occurrence of mislabeling is unsuspected, then $r^*(\mathbf{x})$ will be compared with $k$ rather than $k^*$. In the case of equal misclassification costs, so that $k = 1$, we have $k^* = k = 1$, so that no problems arise from the misclassification. (Indeed, advantages can even arise: see [9].) However, what happens if $k \neq 1$? It is easy to show that $r^*(\mathbf{x}) > r(\mathbf{x})$ whenever $r(\mathbf{x}) < 1$ and that $r^*(\mathbf{x}) < r(\mathbf{x})$ whenever $r(\mathbf{x}) > 1$. That is, the effect of the errors in class labels is to shrink the posterior class odds toward 1, so that comparing $r^*(\mathbf{x})$ with $k$ rather than $k^*$ is likely to lead to worse performance. There is also a secondary issue, that the shrinkage of $r(\mathbf{x})$ will make it less easy to estimate the decision surface accurately because it is a flatter surface: the variance of the estimated decision surface, from sample to sample, will be greater when there is mislabeling of classes. In such circumstances it is better to stick to simpler models, since the higher order terms of the more complicated models will be very inaccurately estimated.

### 4.2 Arbitrariness in the Class Definition

The classical supervised classification paradigm also takes as fundamental the fact that the classes are well defined. That is, that there is some fixed clear external criterion which is used to produce the class labels. In many situations, however, this is not the case. In particular, when the classes are defined by thresholding a continuous variable, then there is always the possibility that the defining threshold might be changed. Once again, this situation arises in consumer credit, where it is common to define a customer as "defaulting" if they fall three months in arrears with repayments. This definition, however, is not a qualitative one (contrast has a tumor/does not have a tumor) but is very much a quantitative one. It is entirely reasonable that alternative definitions (e.g., four months in arrears) might be more useful if economic conditions were to change. This is a simple example, but in many situations much more complex class definitions based on logical combinations of numerical attributes, split at fairly arbitrary thresholds, are used. For example, student grades are often based on levels of performance in continuous assessment and examinations. In detecting vertebral deformities in studies of osteoporosis, the ranges of the anterior, posterior and mid heights of the vertebra, as well as functions of these, such as ratios, are combined in quite complicated Boolean conditions to provide the definition (e.g., [10]). Definitions formed in this sort of way are particularly common in situations that involve customer management. For example, Lewis [31] defined a good account in a revolving credit operation (such as a credit card) as someone whose billing account shows (a) on the books for a minimum of 10 months, (b) activity in 6 of the most recent 10 months, (c) purchases of more than $50 in at least 3 of the past 24 months and (d) not more than once 30 days delinquent in the past 24 months. A bad account is defined as (a) delinquent for 90 days or more at any

time with an outstanding undisputed balance of $50 or more, (b) delinquent three times for 60 days in the past 12 months with an outstanding undisputed balance on each occasion of $50 or more or (c) bankrupt while the account was open. Li and Hand [32] gave an even more complicated example from retail banking.

Our concern with these complicated definitions is that they are fairly arbitrary: the thresholds used to partition the various continua are not natural thresholds, but are imposed by humans. It is entirely possible that, retrospectively, one might decide that other thresholds would have been better. Ideally, under such circumstances, one would go back to the design data, redefine the classes and recompute the classification rule. However, this requires that the raw data have been retained at the level of the underlying continua used in the definitions. This is often not the case. The term *concept drift* is sometimes used to describe changes to the definitions of the classes. See, for example, the special issue of *Machine Learning* (1998, Vol. 32, No. 2), Widmer and Kubat [46] and Lane and Brodley [30]. The problem of changing class definitions has been examined in [25, 26] and [28].

If the very definitions of the classes may change between designing the classification rule and applying it, then clearly there is little point in developing an overrefined model for the class definition which is no longer appropriate. Such models fail to take into account all sources of uncertainty in the problem. Of course, this does not necessarily imply that simple models will yield better classification results: this will depend on the nature of the difference between the design and application class definitions. However, there are similarities to the overfitting issue. Overfitting arises when a complicated model faithfully reflects aspects of the design data to the extent that idiosyncrasies of that data, rather than merely of the distribution from which the data arose, are included in the model. Then simple models, which fit the design data less well, lead to superior classification. Likewise, in the present context, a model optimized on the design data class definition is reflecting idiosyncrasies of the design data which may not occur in application data, not because of random variation, but because of the different definitions of the classes. Thus it is possible that models which fit the design data less well will do better in future classification tasks.

The possibility of arbitrariness in the class definition discussed in this section is quite distinct from the possibility of class priors or relative misclassification costs being changed—referred to in the quote from Provost and Fawcett [36] above—but the possibility of these changes, also, casts doubt on the wisdom of modeling the problem too precisely, that is, of using models which are too sophisticated.

### 4.3 Optimization Criteria and Performance Assessment

When fitting a model to a design set, one optimizes some criterion of goodness of fit (perhaps modified by a penalization term to avoid overfitting) or of classification performance. Many such measures are in use, including likelihood, misclassification rate, cost-weighted misclassification rate, Brier score, log score and area under the ROC curve. Unfortunately, it is not difficult to contrive data sets for which different optimization criteria lead to (e.g.) linear decision surfaces with very different orientations (even to the extent of being orthogonal). Benton [[2], Chap. 4] illustrated this for several real data sets. Clearly, then it is important to specify the criterion to be used when building a classification rule. If the use to which the model will be put is well specified to the extent that a measure of performance can be precisely defined, then this measure should determine the criterion of goodness of fit. All too often, however, there is a mismatch between the criterion used to choose the model, the criterion used to evaluate its performance, and the criterion which actually matters in real application. For example, a common approach might be to use likelihood to estimate a model's parameters, use misclassification rate to assess its performance and use some cost-weighted misclassification rate in practice (e.g., some combination of specificity and sensitivity). In circumstances such as these, it would clearly be pointless to refine the model to a high degree of accuracy from a likelihood perspective, when this may be only weakly related to the real performance objective.

Having said that, one must acknowledge that often precise details of how performance is to be measured in the future cannot be given. For example, in most applications it is difficult to give more than general statements about the relative costs of different kinds of misclassifications. In such cases it might be worthwhile to choose a criterion that is equivalent to averaging over a range of possible costs: likelihood, the area under a receiver operating characteristic curve and the weighted version of the latter described in [1] can all be regarded as attempts to do that.

## 5. INTERPRETING EMPIRICAL COMPARISONS

There have been a great many empirical comparisons of the performance of different kind of classification rules. Some of these are in the context of a new method having been developed and the effort to gain some understanding of how it performs relative to existing methods. Other comparisons are purely comparative studies, seeking to make disinterested comparative statements about the relative merits of different methods. At first glance, such comparative studies are useful in shedding light on the different methods, on which generally yield superior performance or on which are to be preferred for particular kinds of data or in particular domains. However, on closer examination, such comparisons have major weaknesses and can even be seriously misleading. Various authors have drawn attention to these problems, including Duin [4], Salzberg [40], Hand [13], Hoadley [21] and Efron [5], so we will only briefly mention some of the main points here; in particular, only those points relative to classification accuracy, rather than other aspects of performance. Jamain and Hand [24] also gave a more detailed review of comparative studies of classification rules.

Different categories of users might be expected to obtain different rankings of classification methods in comparative studies. For example, we can contrast an expert user, who will be able to fine-tune methods, with an inexperienced user, perhaps someone who has simply pulled some standard public-domain software from the web. It would probably be surprising if their rankings did not differ. Moreover, experts will tend to have particular expertise with particular classes of method. Someone expert in neural networks may well achieve superior results with those methods than with support vector machines and vice versa. Taken to an extreme, of course, many comparative studies are made to establish the performance and properties of newly invented methods—by their inventors. One might expect substantial bias in favor of the new methods, compared to what others might be able to achieve, in such studies. Duin [4] pointed out the difficulty of comparing, "in a fair and objective way," classifiers which require substantial input of expertise (so that domain knowledge can be taken advantage of) and classifiers which can be applied automatically with little external input of expertise. The two extremes (of what is really a continuum, of course) are appropriate in different circumstances.

The principle of comparing methods by applying them to a collection of disparate real data sets is useful, but has its weaknesses. An obvious one is that different studies use different collections of data sets, so making comparisons difficult. Furthermore, the collection will not be representative of real data sets in any formal sense. Moreover, a potential user is not really interested in some "average performance" over distinct types of data, but really wants to know what will be good for his or her problem, and different people have different problems, with data arising from different domains. A given method may be very poor on most kinds of data, but very good for certain problems.

The widespread use of standard collections of data sets (such as the UCI repository [35]) has clear merits: new methods can be compared with earlier ones on a level playing field. However, this also means that there will be some overfitting both to the individual data sets in the collection and to the collection as a whole. That is, some methods will do well on data sets in the collection purely by chance. Indeed, the more successful the collection is in the sense that more and more people use it for comparative assessments, the more serious this problem will become.

Jamain and Hand [24] pointed out the difficulty of saying exactly what a classification "method" is. Is a neural network with a single hidden node to be regarded as from the same family as one with an arbitrary number of hidden nodes? It is clearly not *exactly* the same method. Comparative evaluations using the two models may well yield very different classification results. It is this sort of phenomenon which explains why the comparative performance literature contains many different results for "the same" methods applied to given public data sets. Can one then draw general conclusions about the effectiveness of the method of neural networks? Furthermore, to what extent is preprocessing the data to be regarded as part of the method? Linear discriminant analysis on raw data may yield very different results from the same model applied to data which has been processed to remove skewness. Is, then, linear discriminant analysis good or bad on these data? Likewise, is a data set in which missing values have been replaced by imputed values the same as a data set in which incomplete records have been dropped? Applying the same method to the two variants of the data is likely to yield different results.

We have already commented that the "accuracy" of a classification rule can be measured in a wide variety of ways, and that different measures are likely to yield different performance rankings of classifiers.

Given all of the above points, it is not surprising that different authors have drawn different conclusions about the relative accuracy of different classifiers. Other commentators have taken things even further. In the discussion that accompanies [3], Efron suggested that new methods always look better than older ones and that complicated methods are harder to criticize than simpler ones. He also noted that it is difficult to make fair comparisons by making the same effort in applying different methods—a point made above. Hoadley, in the same discussion, "coined a phrase called the 'ping-pong theorem.' This theorem says that if we revealed to Professor Breiman the performance of our best model and gave him our data, then he could develop an algorithmic model using random forests, which would outperform our model. But if he revealed to us the performance of his model, then we could develop a segmented scorecard, which would outperform his model."

With so many difficulties in ranking and comparing classifiers, one might naturally have reservations about small differences in performance—of the kind generally asserted for the more complicated and sophisticated methods over the older and simpler models.

## 6. CONCLUSION

In Section 2 we demonstrated that, when building predictive models of increasing complexity, the marginal gain from complicated models is typically small compared to the predictive power of the simple models. In many cases, the simple models accounted for over 90% of the predictive power that could be achieved by "the best" model we could find. Now, in the idealized classical supervised classification paradigm, certain assumptions are implicit: it is assumed that the distributions from which the design points and the new points are drawn are the same, that the classes are well defined and the definitions will not change, that the costs of different kinds of misclassification are known accurately, and so on. In real applications, however, these additional assumptions will often not hold. This means that apparent small (laboratory) gains in performance might not be realized in practice—they may well be swamped by uncertainties arising from mismatches between the apparent problem and the real problem. In particular, many of the comparative studies in the literature are based on brief descriptions of data sets, containing no background information at all on such possible additional sources of variation due to breakdown of implicit assumptions of the kind illustrated above. This must cast doubt on the validity of their conclusions. In general, it means that deeper critical assessment of the context of the problem and data should be made if useful practical conclusions are to be drawn. If enough is known about likely additional sources of variability, beyond the classical sources of sampling variability and model uncertainty, then more sophisticated models can be built. However, if insufficient information is known about these additional sources, which we speculate will very often be the case, then the principle of parsimony suggests that it is better to stick to simple models.

We should note, parenthetically, that there are also other reasons to favor simple models. Interpretability, in particular, is often an important requirement of a classification rule. Indeed, sometimes it is even a legal requirement (e.g., in credit scoring). This leads us to the observation that what one regards as "simple" may vary from user to user: some might favor weighted sums of predictor values, others might prefer (small) tree structures and yet others might regard nearest neighbor methods as being simple.

Perhaps it is appropriate to conclude with the comment that, by arguing that simple models are often more appropriate than complex ones and that the claims of superior performance of the more complex models may be misleading, I am not suggesting that no major advances in classification methods will ever be made. Such a claim would be absurd in the face of developments such as the bootstrap and other resampling approaches, which have led to significant advances in classification and other statistical models. All I am saying is that much of the purported advance may well be illusory. Furthermore, although (almost by definition) one cannot predict where the next step-change will come from, one might venture a guess as to its general area. Resampling methods are children of the computer revolution, as indeed are most other recent developments in classifier technology [e.g., classification trees, neural networks, support vector machines, random forests, multivariate adaptive regression splines (MARS) and practical Bayesian methods]. Since progress in computer hardware is continuing, one might reasonably expect that the advances will arise from more powerful data storage and processing ability.

## ACKNOWLEDGMENTS

I have given several presentations based on the ideas in this paper. An earlier version of this paper was presented at the 2004 Conference of the International Federation of Classification Societies in Chicago and appeared in the proceedings [18]. I would like to thank all those who commented on the material. In particular, I am grateful to Svante Wolde, Willi Sauerbrei, Foster Provost, Jerome Friedman and Leo Breiman. Of course, merely because they had valuable and interesting things to say about the ideas does not necessarily mean they agree with them.